# Klein Group And Four Color Theorem

SERGEY KURAPOV

**Abstract.** In this work methods of construction of cubic graphs are analyzed and a theorem of existence of a colored disc traversing each pair of linked edges belonging to an elementary cycle of a planar cubic graph is proved.

## 1. Maximal planar graph

Let *G* be an arbitrary planar connected graph with each vertex having degree of more than two (Fig.1). Let us by adding extra edges split each non-triangular face of the graph *G* into triangular faces (Fig.2). In the result, we obtain a maximal planar graph with all faces being triangles.

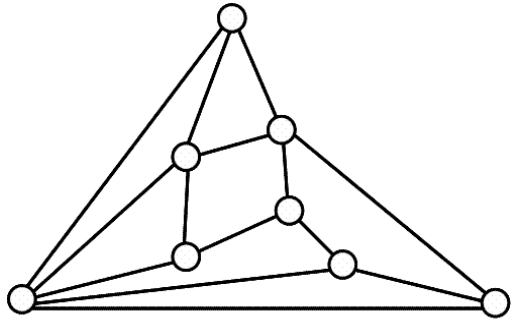 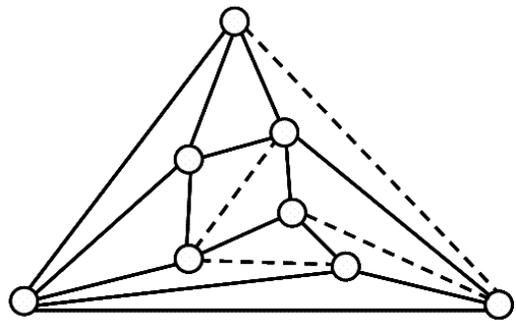

Fig.1. Planar graph *G*.  Fig.2. Maximal planar graph *G*.

Let us inside of each face *s* of the maximal planar graph *G* put a vertex $x^0$, for each edge *u* of the graph *G* let us draw its corresponding edge $u^0$ which connects vertices $x^0_i$ and $x^0_j$ inside faces $s_i$ and $s_j$ on both sides of the edge *u*. Thus constructed topological graph $G^0$ is planar, connected and is named the dual of the graph *G*. The dual graph is an isomorphic graph *H* of degree three or a cubic graph *H*.

Construction of the dual planar cubic graph *H* is shown on Fig.3, where its vertices are colored darker and its edges are shown as dot lines.

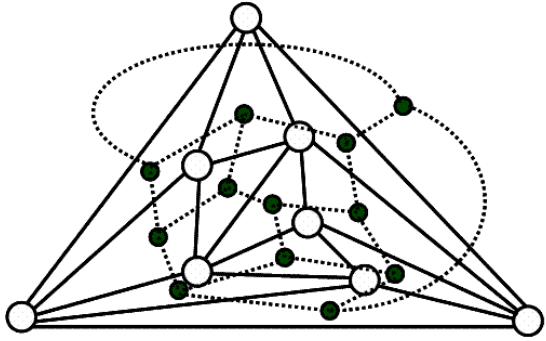
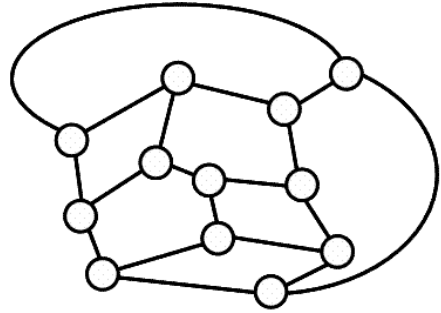

Fig.3. Construction of the planar cubic graph *H*.  Fig.4. Planar cubic graph *H*.

The number of edges in a cubic graph is determined as:

$$m = 3n/2 \qquad (1)$$

Hence, the number of edges in such a graph is always a multiple of three.

Since the number of edges *m* is integer, therefore, the number of vertices *n* in an isomorphic cubic graph is even. Transition from the maximal planar graph *G* with triangular faces to the dual cubic graph *H* is performed by replacing the cycle matrix of the maximal planar graph *G* with the incidence matrix of the cubic graph *H*.

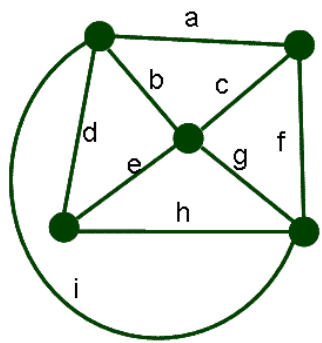
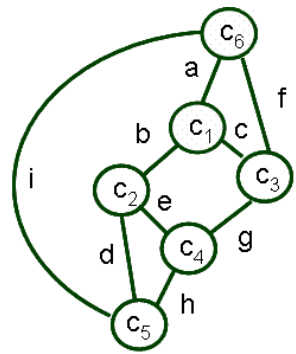

Fig.5. Graph *G*.  Fig.6. Graph *H*.

Cycle matrix of the graph G

|     | a | b | c | d | e | f | g | h | i |
|-----|---|---|---|---|---|---|---|---|---|
| $c_1$ | 1 | 1 | 1 |   |   |   |   |   |   |
| $c_2$ |   | 1 |   | 1 | 1 |   |   |   |   |
| $c_3$ |   |   | 1 |   |   | 1 | 1 |   |   |
| $c_4$ |   |   |   |   | 1 |   | 1 | 1 |   |
| $c_5$ |   |   |   | 1 |   |   |   | 1 | 1 |
| $c_6$ | 1 |   |   |   |   | 1 |   |   | 1 |
|     | a | b | c | d | e | f | g | h | i |

Incidence matrix of the graph H

## 2. Proper coloring of a cubic graph

The following theorem states the coloring relation between the maximal planar graph *G* and its dual planar cubic graph *H*.

**Theorem 1.** (Tait) [1]. A maximal planar graph can be face-colored with 4 colors if and only if its dual graph has chromatic class 3.

If the chromatic class of the cubic graph *H* is 3, then its edges are 3-colorable. Let us define such coloring as the proper coloring of the cubic graph *H*.

The proper coloring implies that edges incident with each vertex of the graph *H* are colored differently. For example, the graph on Fig.7 has proper coloring.

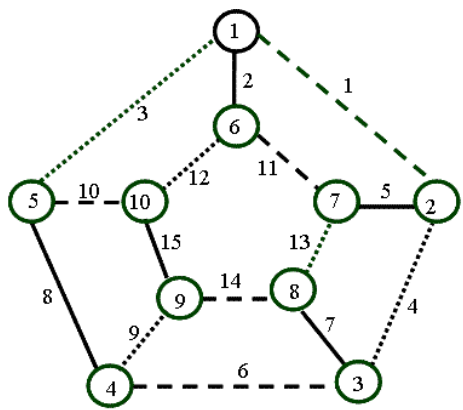

Fig.7. Properly colored cubic graph.

Set of colors in a properly colored cubic graph and the operation of addition yield Klein four-group. The color addition table for this group is as follows:

| + | 0 | 1 | 2 | 3 |
|---|---|---|---|---|
| 0 | 0 | 1 | 2 | 3 |
| 1 | 1 | 0 | 3 | 2 |
| 2 | 2 | 3 | 0 | 1 |
| 3 | 3 | 2 | 1 | 0 |

Here, "0" denotes white colour, let us assign it letter $W$ and use dash-and-dot lines to draw green edges.

"1" denotes red colour, let us assign it letter $R$ and use solid lines to draw red edges.

"2" denotes blue colour, let us assign it letter $B$ and use dot lines to draw blue edges.

"3" denotes green colour, let assign it letter $G$ and use dash lines to draw green edges.

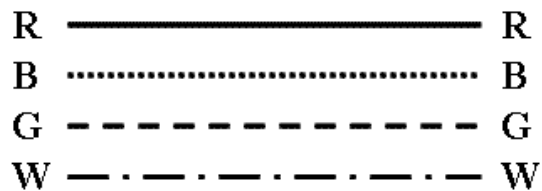

Let us define a quasi-cycle as a cycle with even degree of vertices. Let us define a Hamiltonian quasi-cycle as a cycle with even degree of vertices which crosses each vertex of the graph.

In its turn, a colored Hamiltonian quasi-cycle may consist of a few simple cycles. Let us define each simple cycle as a colored disc. Fig.8d shows two discs, the first traverses edges $\{u_2,u_3,u_8,u_9,u_{12},u_{15}\}$ and the second traverses edges $\{u_4,u_5,u_7,u_{13}\}$.

Let us define an $n$-factor as a regular subgraph of degree $n$. Then:

**Theorem 2.** (Petersen) [3]. Every bridgeless cubic graph is a sum of a 1-factor and a 2-factor.

A properly colored cubic graph always has a color configuration consisting of three colored 1-factors and three colored 2-factors, each 2-factor having discs of even length only.

Every Hamiltonian quasi-cycle (2-factor) consists of a few discs. If there is just one disc, then it is a Hamiltonian cycle.

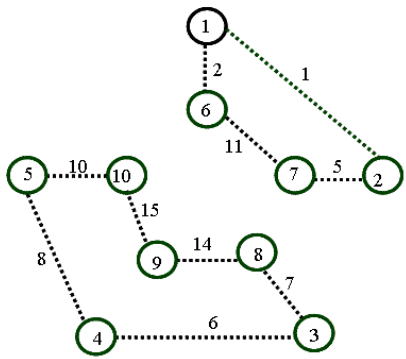

a) blue Hamiltonian quasi-cycle

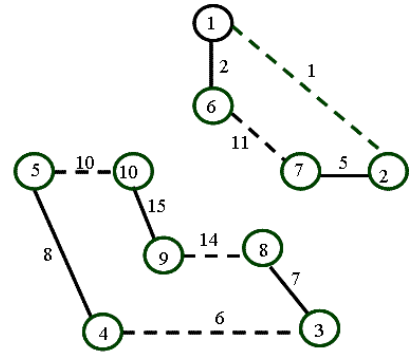

b) the same quasi-cycle consisting of red and green edges

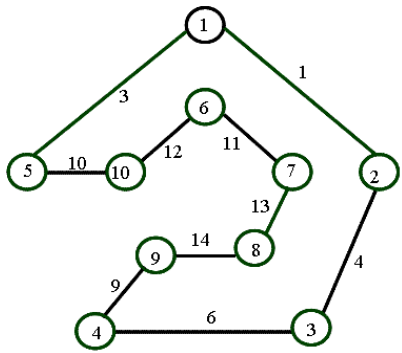

c) red Hamiltonian quasi-cycle

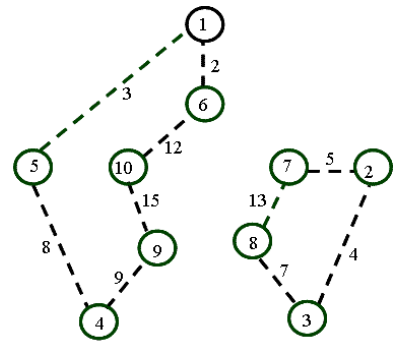

d) green Hamiltonian quasi-cycle

Fig.8. Colored induced Hamiltonian quasi-cycles.

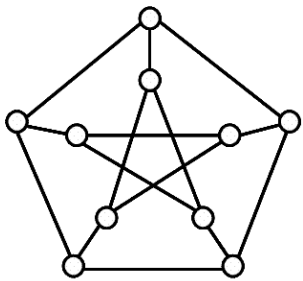

a) The Petersen graph

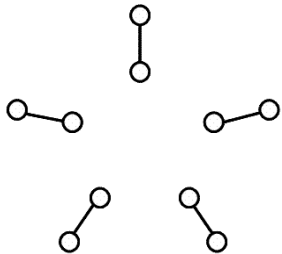

b) The Petersen graph's 1-factor

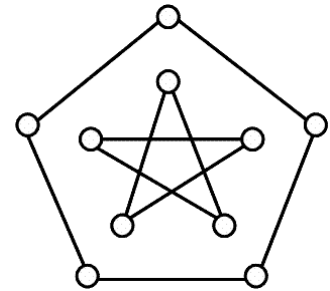

c) The Petersen graph's 2-factor

Fig.9. The Petersen cubic graph.

# 3. Properties of properly colored planar cubic graphs

Colored Hamiltonian quasi-cycles have the following properties:

**Property 1.** Each colored disc has even number of vertices.

Indeed, for a disc to be colored its edges should be 2-colored which is only possible if the number of vertices is even.

**Property 2.** Each colored Hamiltonian quasi-cycle (2-factor) induces two other colored Hamiltonian quasi-cycles (2-factors) and causes proper coloring of the graph $H$.

Given a colored Hamiltonian quasi-cycle, its edges can be 2-colored, uncolored edges can be colored with the color of the Hamiltonian quasi-cycle itself. This makes the cubic graph $H$ colored properly and induces three colored Hamiltonian quasi-cycles.

**Property 3.** Modulo 2 sum of three colored Hamiltonian quasi-cycles (2-factors), induced by the proper coloring of the cubic graph $H$, is an empty set.

As there exist two colored Hamiltonian quasi-cycles traversing each edge, so each edge participates twice in the summation, which in case of modulo 2 sum identically equals to zero. Therefore, the sum of colored Hamiltonian quasi-cycles induced by proper coloring of the cubic graph $H$ is an empty set.

Let $G=(X,U;P)$ be a graph with an enumerated set of edges $U = \{u_1, u_2,...,u_m\}$, and let $£_G$ be a collection of all the spanning subgraphs of the graph $G$. Relative to the addition operation (we define it as the circular sum) the collection $£_G$ yields an Abelian group.

$$(X,U_1;P) \oplus (X,U_2;P) = (X,(U_1 \cup U_2) \setminus (U_1 \cap U_2);P) \tag{3}$$

Indeed, $£_G$ is a known groupoid. Let us to each spanning subgraph of $G = (X,U;P)$ assign a row of numbers $(\alpha_1,\alpha_2,...,\alpha_i,...,\alpha_m)$, where $i = (1,2,...,m)$, and perform row addition as componentwise modulo 2 sum. In the result, we obtain the isomorphic groupoid $£_G$ , with elements being various rows of length $m$ consisting of zero's and one's, which undoubtedly yield an Abelian group.

$$\alpha_i = \begin{cases} 1, if\ u_i \in U; \\ 0, if\ u_i \notin U. \end{cases}$$

Hereinafter, the group $£_G$ will be convenient to consider as a linear space over the field of coefficients $GF(2) =\{0,1\}$, which is also called a spanning subgraph's space of the graph $G$. The dimension of the space *dim* $£_G = m$, since the set of elements *(1,0,...,0),(0,1,...,0),....(0,0,...,1)* provides the basis, with each element being single-edged spanning subgraph.

Regarding cycle subspace, the concept of an elementary cycle in a graph was introduced in [2], where a theorem of existence of the basis in the set of elementary cycles was also proved.

**Definition 1.** An elementary cycle of a graph is a simple cycle such that between its each two non-adjacent vertices do not exist any other paths which are shorter than those belonging to the cycle itself [2].

It must be mentioned that for planar graphs the boundary of a face is determined by an elementary cycle.

In work [2] a MacLane's functional on the cycle subspace's basis of the graph *G* was defined, and it was shown that in case of a planar graph it equals to zero.

$$F(C) = \sum_{i=1}^{m}(S_i - 1)(S_i - 2) = \sum_{i=1}^{m} S_i^2 - 3\sum_{i=1}^{m} S_i + 2m \tag{4}$$

Here $S_i$ is the number of elementary cycles traversing the edge *i*. For non-planar graphs the modulo 2 sum $c_0$ of the basis and the rim of a graph is also empty.

Proceeding from the above, colored Hamiltonian quasi-cycles and colored 1-factors have the following properties:

$$R_c \oplus G_c \oplus B_c = \emptyset; \tag{5}$$

$$R_f \oplus G_f \oplus B_f = H; \tag{6}$$

$$R_f \oplus R_c = H; \tag{7}$$

$$G_f \oplus G_c = H; \tag{8}$$

$$B_f \oplus B_c = H. \tag{9}$$

Here, *c* denotes the set of edges belonging to a 2-factor (Hamiltonian quasi-cycle), *f* denotes the set of edges belonging to a 1-factor, and *H* is the set of edges of the original planar cubic graph.

Futhermore, in a properly colored cubic graph each colored 2-factor may be represented as a finite sum of elementary cycles.

$R_c$ - the sum of elementary cycles inducing a red 2-factor,

$B_c$ - the sum of elementary cycles inducing a blue 2-factor,

$G_c$ - the sum of elementary cycles inducing a green 2-factor,

$W_c$ - the sum of elementary cycles which do not belong to either red, blue or green 2-factors.

$$R_c = \bar{R}_c; \quad B_c = \bar{B}_c; \quad G_c = \bar{G}_c \tag{10}$$

Let us illustrate the above for a properly colored cubic graph $H$ (Fig.7).

The elementary cycles (faces of the planar graph) and the rim of the cubic graph $H$ are:

$c_1 = \{u_1, u_2, u_5, u_{11}\}$, $c_2 = \{u_4, u_5, u_7, u_{13}\}$, $c_3 = \{u_6, u_7, u_9, u_{14}\}$,

$c_4 = \{u_8, u_9, u_{10}, u_{15}\}$, $c_5 = \{u_2, u_3, u_{10}, u_{12}\}$, $c_6 = \{u_{11}, u_{12}, u_{13}, u_{14}, u_{15}\}$,

$c_0 = \{u_1, u_3, u_4, u_6, u_8\}$.

For the proper coloring shown on Fig.7:

$R_c = c_1 \oplus c_2 \oplus c_3 \oplus c_5 = \{u_1, u_2, u_5, u_{11}\} \oplus \{u_4, u_5, u_7, u_{13}\} \oplus \{u_6, u_7, u_9, u_{14}\} \oplus$
$\oplus \{u_2, u_3, u_{10}, u_{12}\} = \{u_1, u_3, u_4, u_6, u_9, u_{10}, u_{11}, u_{12}, u_{13}, u_{15}\}$;

$\bar{R}_c = c_4 \oplus c_6 \oplus c_0 = \{u_8, u_9, u_{10}, u_{15}\} \oplus \{u_{11}, u_{12}, u_{13}, u_{14}, u_{15}\} \oplus \{u_1, u_3, u_4, u_6, u_8\} =$
$= \{u_1, u_3, u_4, u_6, u_9, u_{10}, u_{11}, u_{12}, u_{13}, u_{14}\}$;

$B_c = c_1 \oplus c_3 \oplus c_4 = \{u_1, u_2, u_5, u_{11}\} \oplus \{u_6, u_7, u_9, u_{14}\} \oplus \{u_8, u_9, u_{10}, u_{15}\} =$
$= \{u_1, u_2, u_5, u_6, u_7, u_8, u_{10}, u_{11}, u_{14}, u_{15}\}$;

$\bar{B}_c = c_2 \oplus c_5 \oplus c_6 \oplus c_0 = \{u_4, u_5, u_7, u_{13}\} \oplus \{u_2, u_3, u_{10}, u_{12}\} \oplus \{u_{11}, u_{12}, u_{13}, u_{14}, u_{15}\} \oplus$
$\oplus \{u_1, u_3, u_4, u_6, u_8\} = \{u_1, u_2, u_5, u_6, u_7, u_8, u_{10}, u_{11}, u_{14}, u_{15}\}$;

$G_c = c_2 \oplus c_4 \oplus c_5 = \{u_4, u_5, u_7, u_{13}\} \oplus \{u_8, u_9, u_{10}, u_{15}\} \oplus \{u_2, u_3, u_{10}, u_{12}\} =$
$= \{u_2, u_3, u_4, u_5, u_7, u_8, u_9, u_{12}, u_{13}, u_{15}\}$;

$\bar{G}_c = c_1 \oplus c_3 \oplus c_6 \oplus c_0 = \{u_1, u_2, u_5, u_{11}\} \oplus \{u_6, u_7, u_9, u_{14}\} \oplus \{u_{11}, u_{12}, u_{13}, u_{14}, u_{15}\} \oplus$
$\oplus \{u_1, u_3, u_4, u_6, u_8\} = \{u_2, u_3, u_4, u_5, u_7, u_8, u_9, u_{12}, u_{13}, u_{15}\}$;

$W_c = c_6 \oplus c_0 = \{u_{11}, u_{12}, u_{13}, u_{14}, u_{15}\} \oplus \{u_1, u_3, u_4, u_6, u_8\} = \{u_1, u_3, u_4, u_6, u_8, u_{11}, u_{12}, u_{13}, u_{14}, u_{15}\}$;

$R_c \oplus G_c \oplus B_c = \{u_1, u_3, u_4, u_6, u_9, u_{10}, u_{11}, u_{12}, u_{13}, u_{14}\} \otimes \{u_1, u_2, u_5, u_6, u_7, u_8, u_{10}, u_{11}, u_{14}, u_{15}\} \otimes$
$\oplus \{u_2, u_3, u_4, u_5, u_7, u_8, u_9, u_{12}, u_{13}, u_{15}\} = \emptyset$.

$R_f = \{u_2, u_5, u_7, u_8, u_{15}\}$; $B_f = \{u_3, u_4, u_9, u_{12}, u_{13}\}$; $G_f = \{u_1, u_6, u_{10}, u_{11}, u_{14}\}$;

$R_f \oplus G_f \oplus B_f = H$;

$R_f \oplus R_c = \{u_1, u_3, u_4, u_6, u_9, u_{10}, u_{11}, u_{12}, u_{13}, u_{14}\} \oplus \{u_2, u_5, u_7, u_8, u_{15}\} = H$;

$B_f \oplus B_c = \{u_1, u_2, u_5, u_6, u_7, u_8, u_{10}, u_{11}, u_{14}, u_{15}\} \oplus \{u_3, u_4, u_9, u_{12}, u_{13}\} = H$;

$G_f \oplus G_c = \{u_2, u_3, u_4, u_5, u_7, u_8, u_9, u_{12}, u_{13}, u_{15}\} \oplus \{u_1, u_6, u_{10}, u_{11}, u_{14}\} = H$.

## 4. Rotation of colored discs

Let us define rotation of a colored Hamiltonian disc as a change in coloring order of its edges and denote it by *rot*. Rotation of a disc causes other colored Hamiltonian quasi-cycles to change. Fig.8b shows a blue disc before rotation and Fig.10 shows it after rotation.

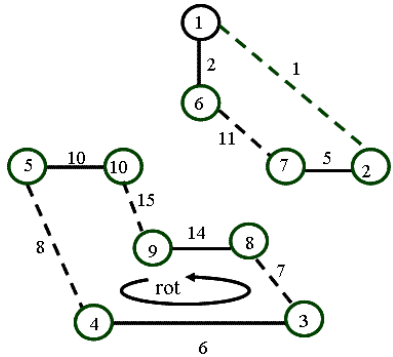

Fig.10. Blue disc after rotation.

The operation of rotation may be described as a permutation of colors of edges in a colored disc. For example, for the graph on Fig.10, the colored 2-factors before rotation are:

blue 2-factor - $\{u_1, u_2, u_5, u_6, u_7, u_8, u_{10}, u_{11}, u_{14}, u_{15}\}$,

red 2-factor - $\{u_1, u_3, u_4, u_6, u_9, u_{10}, u_{11}, u_{12}, u_{13}, u_{14}\}$,

green 2-factor - $\{u_2, u_3, u_4, u_5, u_7, u_8, u_9, u_{12}, u_{13}, u_{15}\}$.

After rotation of the blue disc (Fig.11b):

blue 2-factor - $\{u_1, u_2, u_5, u_6, u_7, u_8, u_{10}, u_{11}, u_{14}, u_{15}\}$,

red 2-factor - $\{u_1, u_3, u_4, u_7, u_8, u_9, u_{11}, u_{12}, u_{13}, u_{15}\}$,

green 2-factor - $\{u_2, u_3, u_4, u_5, u_6, u_9, u_{10}, u_{12}, u_{13}, u_{14}\}$.

In red and green 2-factors permuted pairs of edges are either ($u_6 \leftrightarrow u_7$), ($u_{10} \leftrightarrow u_8$), ($u_{14} \leftrightarrow u_{15}$) or ($u_7 \leftrightarrow u_{14}$), ($u_{10} \leftrightarrow u_{15}$), ($u_8 \leftrightarrow u_6$).

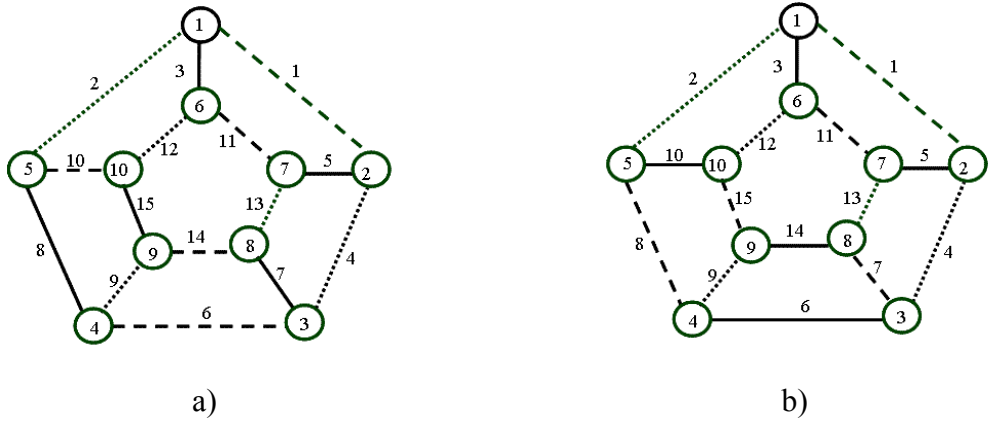

a)  b)

Fig.11 Graph *H* before (a) and after (b) rotation of the blue disc.

## 5. Coloring of elementary cycles

For a properly colored planar cubic graph the circular sum of elementary cycles and the rim is $\sum_{i=1}^{m-n+1} c_i \oplus c_0 = \varnothing$, the circular sum of the three Hamiltonian quasi-cycles is $R_c \oplus G_c \oplus B_c = \varnothing$, therefore, the color of an elementary cycle may be determined from the following equation:

$$(a_{1,1}R + a_{1,2}B + a_{1,3}G + a_{1,4}W)c_1 \oplus (a_{2,1}R + a_{2,2}B + a_{2,3}G + a_{2,4}W)c_2 \oplus$$
$$\oplus \ldots \oplus (a_{m-n+1,1}R + a_{m-n+1,2}B + a_{m-n+1,3}G + a_{m-n+1,4}W)c_{m-n+1} \oplus \qquad (11)$$
$$\oplus (a_{m-n+2,1}R + a_{m-n+2,2}B + a_{m-n+2,3}G + a_{m-n+2,4}W)c_0 ,$$

where $a_{i,j} \in \{0,1\}$; i =(1,2,…,m-n+2); j = (1,2,3,4). Here, coefficients $a_{i,j}$ determine whether an elementary cycle belongs to a colored 2-factor.

Let us illustrate the above for the proper coloring shown on Fig.7.

$$(c_1 \oplus c_2 \oplus c_3 \oplus c_5) \oplus (c_1 \oplus c_3 \oplus c_4) \oplus (c_2 \oplus c_4 \oplus c_5) =$$
$$= (R+B)c_1 \oplus (R+G)c_2 \oplus (R+B)c_3 \oplus (B+G)c_4 \oplus (R+G)c_5 =$$
$$= Gc_1 \oplus Bc_2 \oplus Gc_3 \oplus Rc_4 \oplus Bc_5 .$$

Also, there exists a transition from one to another coloring of elementary cycles. It is performed by adding the chosen color to all colors of the basis of elementary cycles and to the rim of a planar cubic graph.

For example, addition of the *R*-color to the colors of elementary cycles results in the subsequent re-coloring of elementary cycles:

$$(G+R)c_1 \oplus (B+R)c_2 \oplus (G+R)c_3 + (R+R)c_4 + (B+R)c_5 + (W+R)c_6 + (W+R)c_0 =$$
$$Bc_1 + Gc_2 + Bc_3 + Wc_4 + Gc_5 + Rc_6 + Rc_0.$$

And so on. The said is represented by the following table:

|  | $c_1$ | $c_2$ | $c_3$ | $c_4$ | $c_5$ | $c_6$ | $c_0$ |
|---|---|---|---|---|---|---|---|
| Original coloring | G | B | G | R | B | W | W |
| +R | B | G | B | W | G | R | R |
| +B | R | W | R | G | W | B | B |
| +G | W | R | W | B | R | G | G |

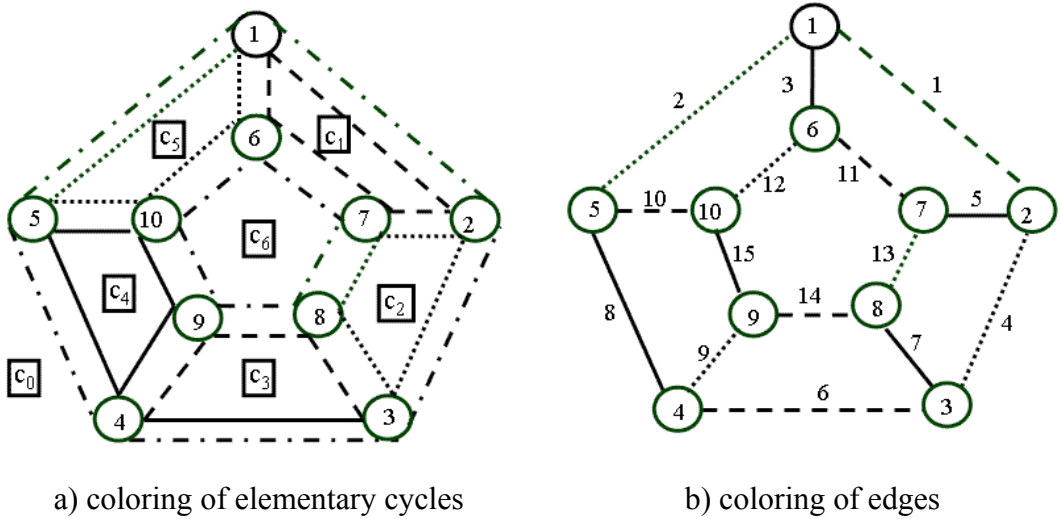

a) coloring of elementary cycles    b) coloring of edges

Fig.12. Coloring of elementary cycles and edges of a planar cubic graph.

Let us denote the sum of red elementary cycles as $\Sigma c_R$, the sum of blue elementary cycles as $\Sigma c_B$, the sum of green elementary cycles as $\Sigma c_G$, the sum of white elementary cycles as $\Sigma c_W$ Then:

$$\Sigma c_R + \Sigma c_B = \Sigma c_G + \Sigma c_W ; \tag{12}$$

$$\Sigma c_R + \Sigma c_G = \Sigma c_B + \Sigma c_W ; \tag{13}$$

$$\Sigma c_G + \Sigma c_B = \Sigma c_R + \Sigma c_W . \tag{14}$$

Equations (11) show that in a planar cubic graph coloring of edges, caused by group addition of two colored Hamiltonian quasi-cycles traversing an edge, on one hand, yields coloring of elementary cycles. On the other hand, the proper coloring of elementary cycles induces colored Hamiltonian quasi-cycles, while group addition of colors of two elementary cycles traversing an edge, yields its coloring.

## 6. Principles of construction of cubic graphs

Every cubic graph on the plane can me constructed by adding new edges using the following methods:

**Method 1.** Put two new vertices on two edges of the original cubic graph and draw a new edge which connects these vertices (Fig.13).

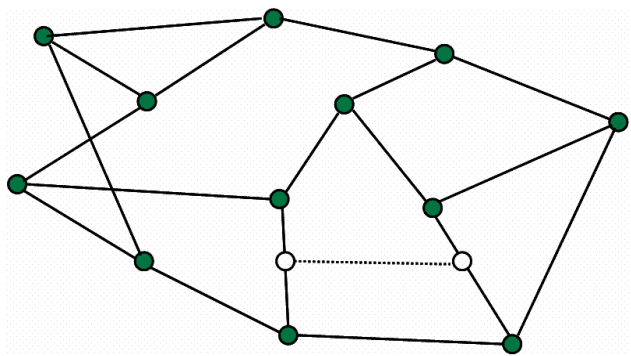

Fig.13. Construction of the cubic graph *H* (method 1).

**Method 2.** Put two new vertices on one edge of the original cubic graph and draw a new edge which connects these vertices (Fig.14).

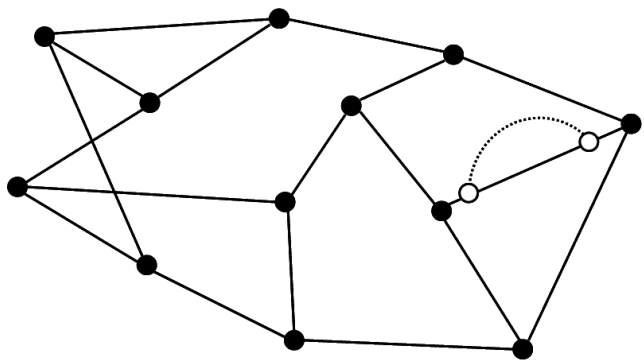

Fig.14. Construction of the cubic graph *H* (method 2).

The minimal graph necessary for the initial construction is the cubic graph with 2 vertices and 3 edges (Fig.15).

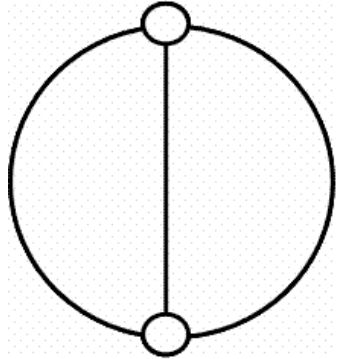

Fig.15.Minimal cubic graph.

Each new cubic graph may be constructed from the preceding cubic graph by adding a new edge using one of the two methods. We may notice that the color of the new edge conforms to the color of the Hamiltonian disc the two vertices of the new edge belong to. Indeed, addition of the new edge increases the number of vertices exactly by 2 which, of course, does not affect the coloring of the graph. This also holds for removing of edges, as it decreases the number of vertices in a disc exactly by 2. Therefore, in case of adding or removing an edge, the proper coloring of a cubic graph is feasible.

Let us define a planar cubic graph before adding a new edge as *preceding*, a planar cubic graph after the new edge has been added as *resulting*. Edges of the preceding graph which own vertices of the new edge let us define as *linked edges*. In order to obtain the resulting planar cubic graph linked edges should belong to an elementary cycle.

Then coloring of any planar cubic graph may be reduced to the following sequential procedure:
- removing edges from the original planar cubic graph to obtain the planar cubic graph with a known coloring;
- coloring of the preceding planar cubic graph;
- searching for the linked edges to add a new edge;
- searching for the colored disc which traverses the linked edges;
- coloring of the new edge using the color of the disc which traverses the linked edges;
- re-coloring of edges of the disc which traverses the linked edges after the new edge has been added.

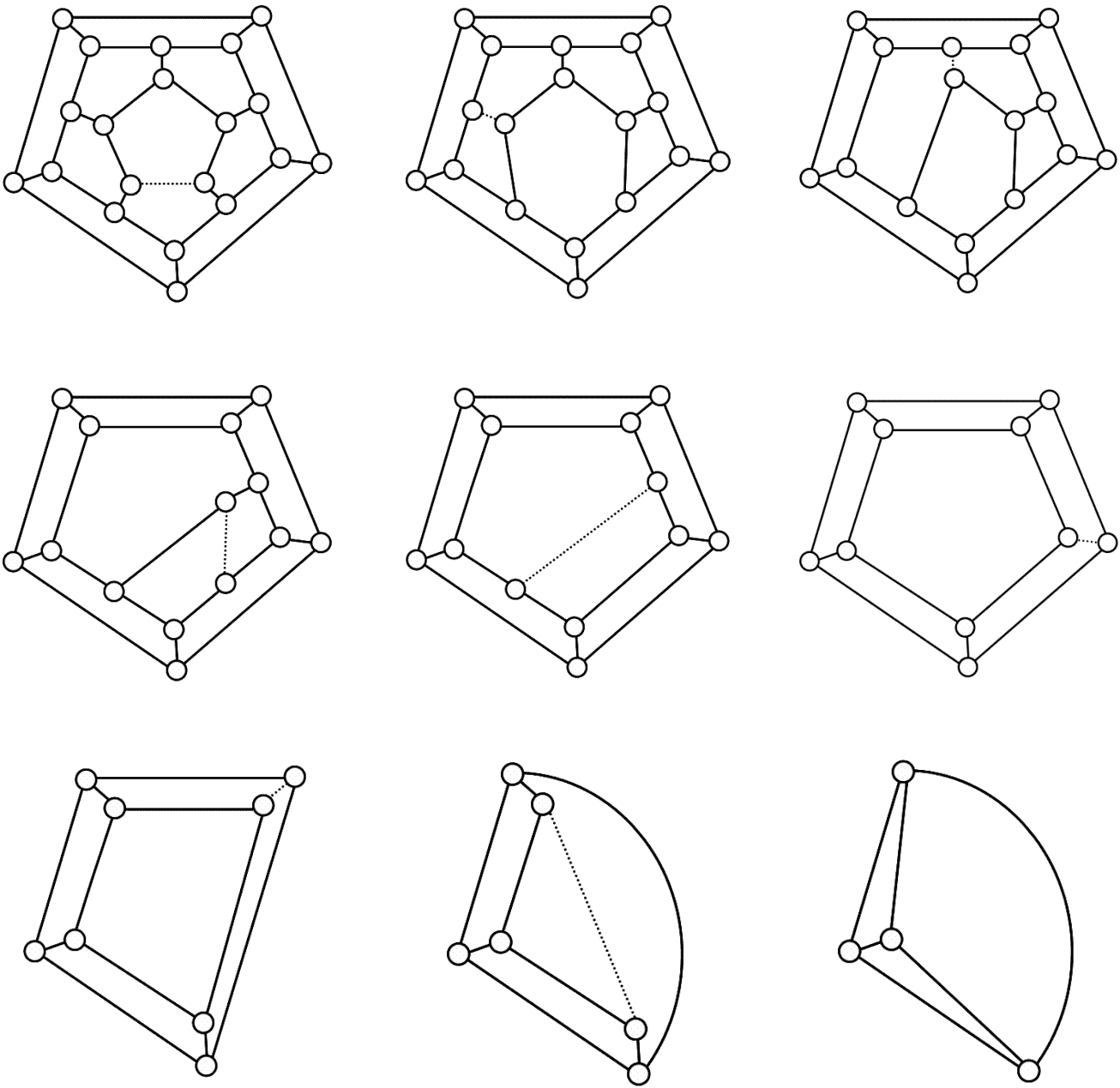

Fig 16. Step-by-step procedure of removing edges from the original planar cubic graph.

Fig.16 shows the step-by-step procedure of removing edges from the original cubic graph with an aim to obtain a known properly colored planar cubic graph (edges being removed are shown as dot lines)

Fig.17 shows the step-by-step procedure of adding new edges and subsequent proper coloring with an aim to obtain the original properly colored planar cubic graph from a known properly colored planar cubic graph.

As there are two different colored discs traversing each linked edge, the following combinations are possible:

1. both linked edges belong to a common colored disc;
2. both linked edges belong to two different colored discs;
3. each linked edge belongs to two different colored discs of one color and two different colored discs of another color;
4. each linked edge belongs to two different discs of the same first color, one linked edge belongs to a disc of the second color and another linked edge belongs to the disc of the third color.

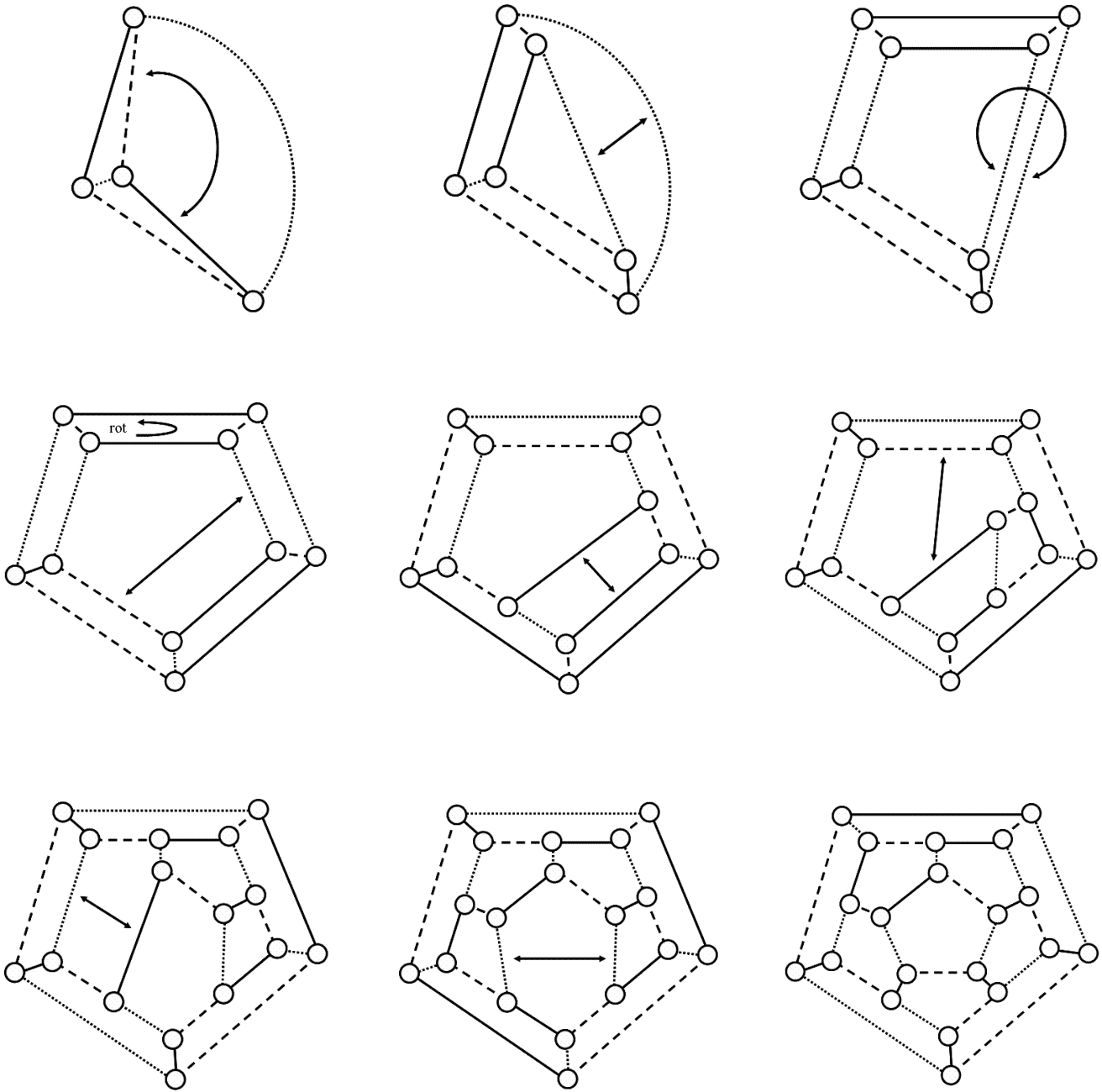

Fig.17 Step-by-step procedure of adding new edges and proper coloring.
(the double-arrow line ↔ represents the new edge being added)

It is easy to show that case 4 is easily reduced to case 3 or case 1 by rotation of the discs of the first color (Fig.18). Thus, in cases 1 and 2 the new edge can be easily colored. Let us prove, that it may be colored in case 3 too.

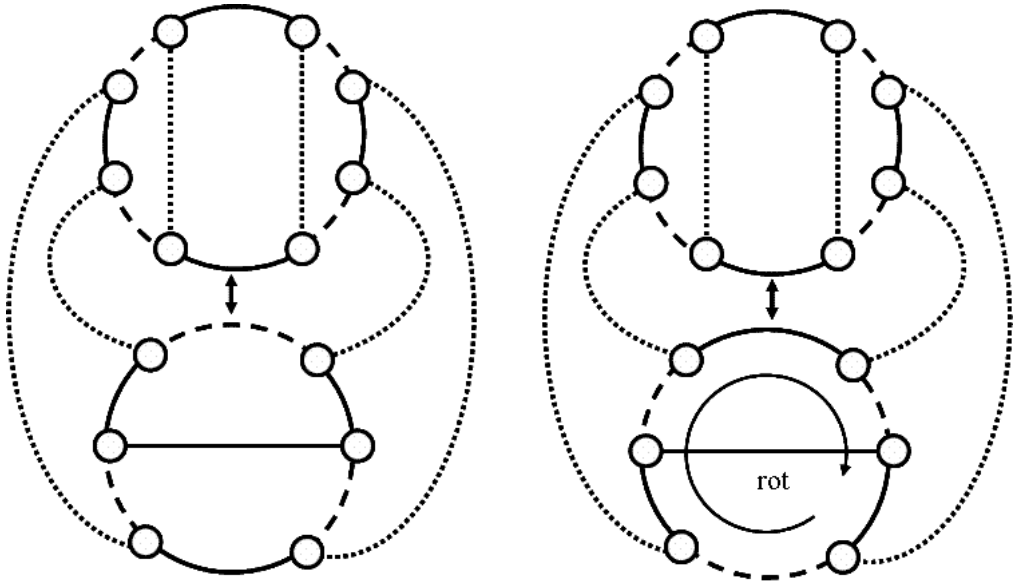

Fig.18. Reduction of case 4 to case 3.
(the double-arrow line ↔ represents the new edge being added)

## 7. Theorem of coloring of a newly added edge.

Suppose that a colored disc traversing the pair of linked edges doesn't exist. Then after adding of a new edge, half of each linked edge has to be colored white to satisfy coloring rules (for every vertex incident edges should be colored differently). In this case, every disc traversing the linked edges increases in length by 1 which makes its length odd (Fig.19). But then there should exist three differently colored elementary discs traversing white-colored edges to satisfy the group addition rule:

$$W = R + B + G.$$

But the latter assertion breaks the rule of planarity (according to MacLane's theorem, the graph is planar if and only if its cycle space has a simple basis). We have come to a contradiction. Hence, a planar cubic graph always has a colored disc traversing the linked edges which belong to a common elementary cycle.

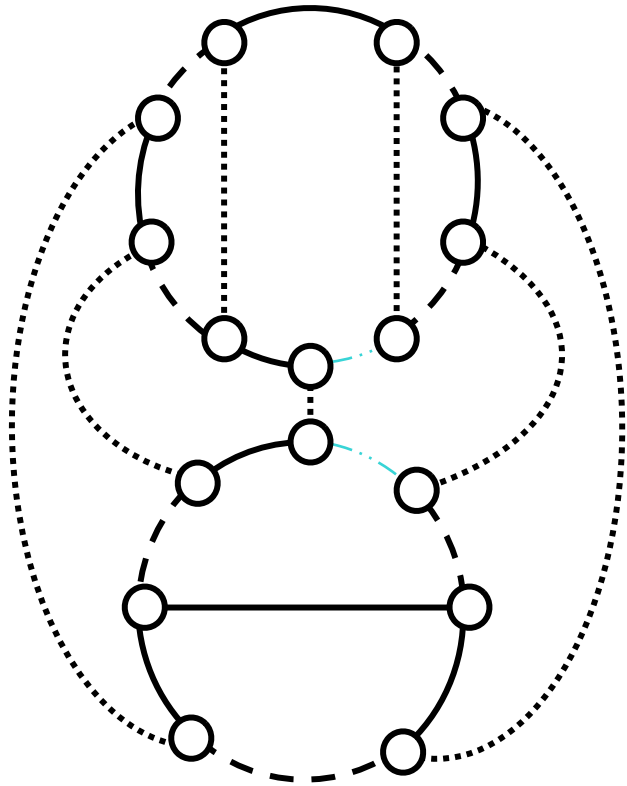

Fig.19. Adding white edges to the planar graph.

Thus, the folllowing theorem has been proved.

**Theorem.** In a bridgeless planar cubic graph there always exists a colored disc traversing each pair of linked edges belonging to the elementary cycle of the graph.

**Corollary 1.** The newly added edge for a resulting bridgeless planar cubic graph is always colorable.

This follows from the fact that a new edge may always be colored using the color of the disc which traverses the pair of linked edges.

**Corollary 2.** A bridgeless planar cubic graph has chromatic class 3.

It follows from the sequence of the coloring procedure, as coloring of the resulting graph is performed by adding of the new edge and its coloring in the preceding graph with consequent re-coloring of edges of the disc which traverses the linked edges.

**Corollary 3.** In a bridgeless planar cubic graph there always exists a Hamiltonian quasi-cycle (2-factor) with discs of even length.

This follows from the principle of the coloring procedure.

**Corollary 4.** A maximal planar graph has chromatic number 4.

This follows from Corollary 2 and Theorem 1 (Tait).

Let us illustrate our conclusions for a properly colored non-planar graph (Fig.20)

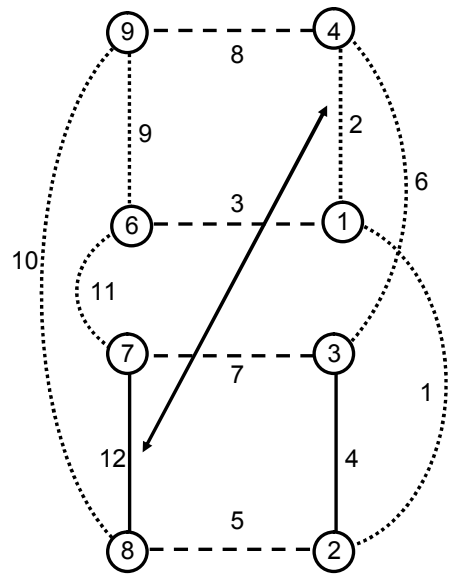

Fig.20. Properly colored non-planar graph.
(the double-arrow line ↔ represents the new edge being added)

For the pair of linked edges ($u_2, u_{15}$) let us add a new edge to obtain the resulting cubic graph. In this non-planar graph there is no colored disc traversing this pair of linked edges. Therefore, in a disc with odd length let us color the edges ($u_2, u_{15}$) of the new graph using white color (Fig.21).

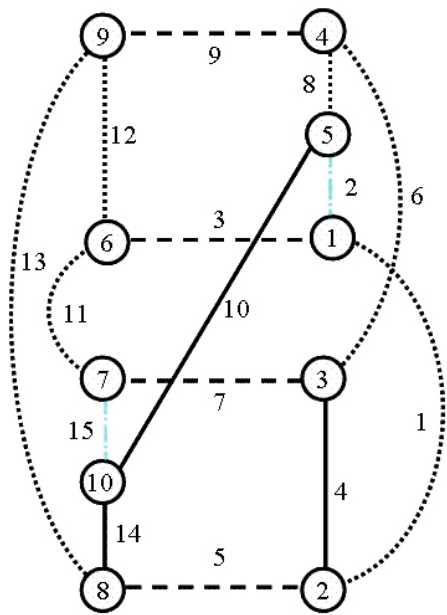

Fig.21. Coloring of edges using white color.
(the double-arrow line ↔ represents the new edge being added)

The set of elementary cycles for this graph consists of twelve cycles with length of 5:

$c_1 = \{u_1,u_2,u_4,u_6,u_8\}$, $c_2 = \{u_1,u_2,u_5,u_{10},u_{14}\}$, $c_3 = \{u_1,u_3,u_4,u_7,u_{11}\}$, $c_4 = \{u_1,u_3,u_5,u_{12},u_{13}\}$,

$c_5 = \{u_2,u_3,u_{10},u_{11},u_{15}\}$, $c_6 = \{u_2,u_3,u_8,u_9,u_{12}\}$, $c_7 = \{u_4,u_5,u_6,u_9,u_{13}\}$, $c_8 = \{u_4,u_5,u_7,u_{14},u_{15}\}$,

$c_9 = \{u_6,u_7,u_9,u_{11},u_{12}\}$, $c_{10} = \{u_6,u_7,u_8,u_{10},u_{15}\}$, $c_{11} = \{u_8,u_9,u_{10},u_{13},u_{14}\}$, $c_{12} = \{u_{11},u_{12},u_{13},u_{14},u_{15}\}$.

Choosing the following cycles as a basis

$c_1 = \{u_1,u_2,u_4,u_6,u_8\}$, $c_3 = \{u_1,u_3,u_4,u_7,u_{11}\}$, $c_5 = \{u_2,u_3,u_{10},u_{11},u_{15}\}$, $c_6 = \{u_2,u_3,u_8,u_9,u_{12}\}$,

$c_8 = \{u_4,u_5,u_7,u_{14},u_{15}\}$, $c_{12} = \{u_{11},u_{12},u_{13},u_{14},u_{15}\}$

indeed proves that there are 3 elementary cycles traversing the linked edges (Fig.22).

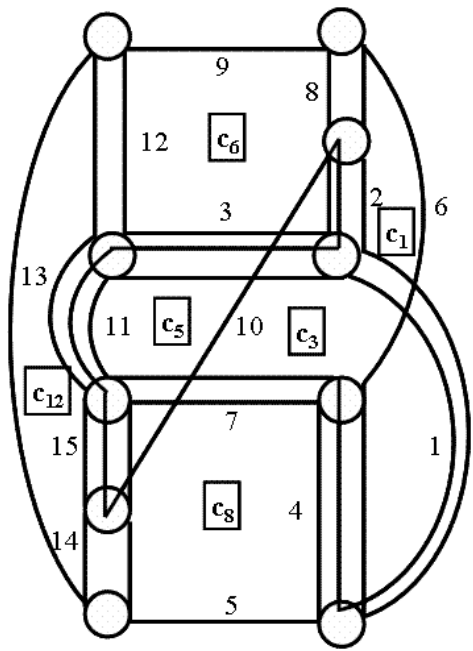

Fig.22. Basis elementary cycles of a non-planar graph.

After re-drawing this graph we obtain a better known representation of it by Petersen (Fig.23).

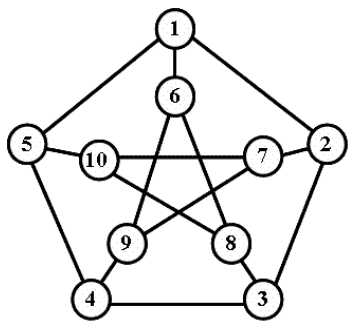

Fig.23. The Petersen graph.

Department of Applied Mathematics
Zaporizhzhya National University
Zaporizhzhya, Ukraine
   *E-mail address:* `kurapov@zsu.zp.ua`